\documentclass[12pt,a4paper]{amsart}
\usepackage{amssymb,graphicx}

\newtheorem{teor}{Theorem}[section]

\newtheorem{cor}[teor]{Corollary}
\newtheorem{prop}[teor]{Proposition}

\newtheorem{ex}[teor]{Example}

    \newcommand{\C}{\mathbb C}
    \newcommand{\R}{\mathbb R}

\begin{document}

\title[Inner Metric Geometry of Complex
Algebraic Surfaces ] {Inner Metric Geometry of Complex Algebraic
Surfaces with Isolated Singularities}

\author[L. Birbrair]{Lev Birbrair}
\author[A. Fernandes]{Alexandre Fernandes}
\address{Departamento de Matem\'atica, Universidade Federal do Cear\'a, Av.
Mister Hull s/n,Campus do PICI, Bloco 914,CEP: 60.455-760 -
Fortaleza - CE - Brasil.} \email{birb@ufc.br}
\email{alex@mat.ufc.br}
\date{\today}

\markboth {L. Birbrair A. Fernandes }{Inner Metric Geometry of
Complex Algebraic Surfaces with Isolated Singularities}



\begin{abstract}
We produce examples of complex algebraic surfaces with isolated
singularities such that these singularities are not metrically
conic, i.e. the germs of the surfaces near singular points are not
bi-Lipschitz equivalent, with respect to the inner metric, to cones.
The technique used to prove the nonexistence of the metric conic
structure is related to a development of Metric Homology. The class
of the examples is rather large and it includes some surfaces of
Brieskorn.
\end{abstract}

\maketitle

\section{Introduction}

An extremely important corollary of the "Triangulation Theorem" of
Lojasiewicz \cite{L} is the existence of topologically conic
structure near a singular point of an algebraic set (real or
complex). Namely, Lojasiewicz proved the following: let $X$ be an
algebraic (or semialgebraic) set in $\R^n$ and let $x_0\in X$. Then
there exists a number $\epsilon >0$ such that the intersection of
$X$ with a ball centered at $x_0$ and radius $\epsilon$ is
homeomorphic to a cone over the intersection of $X$ with the sphere
of the same radius and centered at $x_0$ (the intersection with a
small sphere is usually called the link of $X$ at $x_0$). Moreover,
he proved that the homeomorphism can be chosen as a semialgebraic
map. The same result, but without the conclusion on the
semialgebraicity of the corresponding homeomorphism was obtained by
Milnor \cite{M}, using the integration of the radial vector field,
for complex algebraic hypersurfaces.

One can ask the following. Is the same result true in the sense of
Metric Geometry? Namely, is the intersection of $X$ with a ball of
small radius centered at $x_0$ bi-Lipschitz homeomorphic with
respect to the inner metric to a cone over the intersection of $X$
with a small sphere with the same center?

For real algebraic sets the answer is negative. One can consider
so-called $\beta$-horn, i.e. the algebraic set defined by
$$\{(x_1,x_2,x_3)\in\R^3 \ : \ (x_1^2+x_2^2)^q=x_3^{2p}\}$$
where $p>q$ are coprime positive integers and
$\displaystyle\beta=\frac{p}{q}$. In this case, the corresponding
tangent cone has the real dimension $1$, but the algebraic set
itself has the real dimension $2$. For complex algebraic sets, this
phenomena does not exists, i.e. a tangent cone of complex algebraic
set is a complex algebraic set and has the same dimension of the
corresponding algebraic set \cite{W}.

The question of existence of the inner Lipschitz conic structure is
also very important in the Intersection Homology Theory and
$L_p$-cohomology. If all singularities of a complex algebraic set
$X$ satisfy this property, then, using the results of Brasselet,
Goresky, MacPherson \cite{BGM} and Youssin \cite{Y}, one can show
that these two cohomology theories are isomorphic.

The main goal of this paper is the following. We show that there
exists a big class of complex algebraic surfaces, with isolated
singularities, such that their singular points do not admit a
metrically conic structure. Namely, for any complex algebraic
surface $X$ from this class and for any $\epsilon >0$, the
intersection of $X$ with a ball of the radius $\epsilon$ and
centered at $x_0$ is not bi-Lipschitz equivalent to a cone over any
Nash manifold  (in particular over the intersection with the
corresponding sphere).

We consider two possible versions of the question of the existence
of "Lipschitz conic structure". The first version of the question is
the following: given $X$ an algebraic set in $\R^n$ and $x_0$ a
singular point of $X$, is the intersection of $X$ with a small ball
centered at $x_0$ semialgebraically bi-Lipschitz homeomorphic to the
cone over the intersection of $X$ with a small sphere centered at
$x_0$? When the answer is positive, we say that $X$ is
\emph{strongly metrically conic} at $x_0$. For weight homogeneous
(not homogenous) surfaces in $\C^3$, defined by real polynomials, we
present a criteria on nonexistence of this structure. In fact, we
proved that if the real part of $X$ has empty intersection with the
union of the coordinates hyperplane in $\C^3-\{0\}$, then $X$ is not
strongly metrically conic at $0$. In order to show this result, we
use the theory of "Characteristic Exponents" developed in \cite{BC}
and "Metric Homology Theory" developed in \cite{BB1}, \cite{BB2}. If
the intersection of real part of $X$ with the link of $X$ at $x_0$
presents a nontrivial cycle in $1$-dimensional homology of this
link, we use the methods developed in \cite{BF} to compute a
characteristic exponent of this singularity. The answer, obtained
here, is different then the corresponding answer for strongly
metrically conic singularity (see \cite{BB2}). If the intersection
of real part of $X$ with the link of $X$ at $x_0$ presents a trivial
cycle in $1$-dimensional homology of this link, then we create a
so-called Cheeger's cycle in $3$-dimensional Local Metric Homology
of $X$. According to \cite{BB2}, $3$-dimensional Local Metric
Homology of strongly metrically conic singularity is generate only
by the fundamental cycle of the link. We show that the Cheeger's
cycle and fundamental cycle are independent. Note that the existence
of Cheeger's cycles shows the Filtration Theorem proved in
\cite{BB1} for $1$-dimensional Local Metric Homology is not true for
$3$-dimensional Local Metric Homology. This fact is important for
the further development Metric Homology Theory.

Finally, we consider the second version of the question about
"Lipschitz conic structure". In this case, we do not suppose that
the corresponding homeomorphism is semialgebraic. The singularities
of this sort are called \emph{metrically conic} and not strongly
metrically conic. We present a criteria of nonexistence of this
structure for weighted homogeneous (no homogeneous) surfaces in
$\C^3$ defined by real polynomials. Namely, the image of the real
part of the weighted homogeneous surface by the projection of the
corresponding Seifert Fibration has to have more than one connected
component. We show that the vanishing rate of the "real cycles" is
bigger than one, and we also show that this cannot happen in the
conic case.

In the section 6 we show that there exists  series of the surfaces
of Brieskorn such that the singularity at zero of these surfaces is
also not strongly metrically conic. We show that these surfaces have
Cheeger's cycles. Note that all the Brieskorn surfaces do not
satisfies the conditions of Theorem \ref {theor1} and the conditions
of Theorem \ref {theor2}.

Note that the results of this paper have the following application
in the Theory of Minimal Surfaces. The results of L. Caffarelli, R.
Hardt and L. Simon \cite{CHS} produce the examples of embedded
minimal hypersurfaces in $\R^{n+1}$ ($n\geqslant 3$) with isolated
singularities which are not conic in a direct sense. According to
Federer \cite{F}, compact parts of complex algebraic sets are
area-minimizing, hence our examples are also examples of
area-minimizing with isolated singularities which are not conic even
in a metric sense.

\bigskip

\noindent{\bf Acknowledgements.} We are grateful for Professor W.
Neumann for the extremely important comments in the beginning of
these investigations. We are grateful for Professor T. Mostowski for
the clarification of two different viewpoints on the problems of the
existence of the conic structure and for Professor J. Lira for his
comments on the conic structure of minimal surfaces with
singularities. The first author was supported by CNPq grant N
300985/93-2. The second author was supported by CNPq grant N
300393/2005-9.

\section{Definitions and Notations}
Let $X\subset\R^n$ be a connected algebraic set. We define an
\emph{inner metric} on $X$ as follows. Let $x,y\in X$. The inner
distance $d(x,y)$ is defined as infimum of lengths of rectificable
arcs $\gamma\colon [0,1]\rightarrow X$ such that $\gamma(0)=x$ and
$\gamma(1)=y.$ Note that for connected algebraic sets the inner
metric is well defined.

Let $X\subset\C^n$ be an algebraic set and let $x_0\in X$ be a
singular point. We say that $x_0$ is a \emph{strongly metrically
conic} singular point if for sufficiently small $\epsilon>0$, there
exists a semialgebraic bi-Lipschitz homeomorphism, with respect the
inner metric, $$h:c[X\cap S_{\epsilon}(x_0])\rightarrow X\cap
B_{\epsilon}[x_0],$$ where $B_{\epsilon}[x_0]$ is the closed ball
with center $x_0$ and the radius $\epsilon$, $S_{\epsilon}(x_0)$ is
the sphere with center $x_0$ and the radius $\epsilon$, $c[X\cap
S_{\epsilon}(x_0)]$ is a cone over $X\cap S_{\epsilon}(x_0)$.

The above definition of strongly metrically conic singular point is
equivalent to the following. Consider a semialgebraic triangulation
of $X$ and consider the stars of $x_0,$ according to this
triangulation. The point $x_0$ is strongly metrically conic if the
intersection $X\cap B_{\epsilon}[x_0]$ is semialgebraically
bi-Lipschitz homeomorphic to the star of $x_0$, considered with the
standard metric of the simplicial complex (see \cite{BGM}). An
isolated singular point $x_0\in X$ is called \emph{metrically conic}
if for sufficiently small $\epsilon >0$ there exist a Nash manifold
$N$ and a bi-Lipschitz homeomorphism $$h\colon c[N]\rightarrow X\cap
B_{\epsilon}[x_0],$$ where $c[N]$ is a cone over $N$. Note that we
do not suppose $N$ to be homeomorphic to $X\cap S_{\epsilon}(x_0).$

\section{Weighted Homogeneous Surfaces}

Let $w_1,w_2,w_3$ be positive integer numbers and $w=(w_1,w_2,w_3)$.
Let
$$\alpha_w\colon\C^*\times\C^3\rightarrow\C^3$$ be defined by
$$\alpha_w(t,x)=(t^{w_1}x_1,t^{w_2}x_2,t^{w_3}x_3);$$
where $x=(x_1,x_2,x_3)$. We say that $X\subset\C^3$ is
\emph{weighted homogeneous} with respect to $w=(w_1,w_2,w_3)$ if $X$
is invariant by $\C^*$-action $\alpha_w$. When $w_1=w_2=w_3$, we say
that $X$ is \emph{homogeneous}.

A 2-dimensional complex algebraic subset $X\subset\C^3$ which is
weight homogeneous with respect to $w=(w_1,w_2,w_3)$ is called a
\emph{weight homogeneous algebraic surface in} $\C^3$ with respect
to $w=(w_1,w_2,w_3)$.

\begin{ex} Let $f(X_1,X_2,X_3)\in\C[X_1,X_2,X_3]$ be a nonzero polynomial
and let $w_1,w_2,w_3,d$ be positive integer numbers such that:
$$f(t^{w_1}x_1,t^{w_2}x_2,t^{w_3}x_3)=t^df(x_1,x_2,x_3),$$
$\forall t\in\C^{*}$ and $\forall (x_1,x_2,x_3)\in\C^3.$ Then
$$X=\{(x_1,x_2,x_3)\in\C^3 \ : \ f(x_1,x_2,x_3)=0\}$$ is a weight
homogeneous algebraic surface in $\C^3$ with respect to
$w=(w_1,w_2,w_3)$. In this case, we say that $X$ is defined by a
weight homogeneous polynomial with respect to $w=(w_1,w_2,w_3)$.
\end{ex}

\medskip

In the following, let $X$ be a weight homogeneous algebraic surface
in $\C^3$ with respect to $w=(w_1,w_2,w_3)$, where $w_1,w_2,w_3$ are
coprime positive integers. Let $\varphi\colon\C^3\rightarrow\C^3$ be
defined by
$$\varphi(x_1,x_2,x_3)=(x_1^{w_1},x_2^{w_2},x_3^{w_3})$$ and let
$X^{\prime}={\varphi}^{-1}(X)$. We call $(\varphi,X^{\prime})$
\emph{the homogeneous modification of} $X$.

\medskip

The homogeneous modification $X^{\prime}$ of $X$ is a homogeneous
algebraic surface in $\C^3$, thus it defines a projective complex
algebraic curve
$$M^{\prime}=\{[z_1:z_2:z_3]\in\C P^2 \ : (z_1,z_2,z_3)\in
X^{\prime}\}.$$

Let us consider that $0$ is an isolated singular point of $X$. So,
Let $M$ be defined by $M=(X-\{0\})/{\alpha_w}$ and let $\pi\colon
X\rightarrow M$ be the standard projection. $M$ is a $2$-dimensional
orbifold (see \cite{N}) and $\pi\colon X-\{0\}\rightarrow M$ is a
Seifert Fibration in the sense of Orlik and Wagreich (see
\cite{OW}). Moreover, there exist a branched covering $\phi\colon
M^{\prime}\rightarrow M$ such that
$$\pi\circ\varphi(z_1,z_2,z_3)=\phi[z_1:z_2:z_3]$$ for all
$(z_1,z_2,z_3)\in X^{\prime}-\{0\}$. $(\pi,M)$ is called
\emph{Seifert Fibration} of $X-\{0\}$ associated to
$w=(w_1,w_2,w_3)$ (see \cite{OW}).

\begin{prop}\label{prop_A} Let $F=\{[z_1:z_2:z_3]\in M^{\prime} \ : \
z_1z_2z_3=0\}.$ Then, $M-\phi(F)$ admits a holomorphic structure
such that $$\phi\colon M^{\prime}-F\rightarrow M-\phi(F)$$ is
locally biholomorphic.
\end{prop}

\begin{proof} See \cite{OW}.
\end{proof}

\section{Cheeger's Cycles}

\begin{teor}\label{cheeger's_cycles}
Let $X\subset\R^n$ be a $k$-dimensional semialgebraic set and
$x_0\in X$ be an isolated singular point of $X$ with a connected
local link. Let $Y\subset X$ be a semialgebraic subset satisfying:
\begin{enumerate}
\item $x_0\in Y$ and $X-Y$ has exactly two connected
components $W_1$ and $W_2$;
\item $\mu(X,x_0)=\mu(W_1,x_0)=\mu(W_2,x_0)=k$;
\item $\mu(Y,x_0)>k$.
\end{enumerate}
Then, for $k<\nu <\mu(Y,x_0)$ the space $MH_{loc,k-1}^{\nu}(X,x_0)$
has a subspace isomorphic to $\R^2$.
\end{teor}

\noindent{\bf Remark.} According to the results of \cite{BB2}, the
existence of such a subset $Y\subset X$ proves that the singularity
of $X$ at $x_0$ is not strongly metrically conic.

\medskip

\begin{proof}[Proof of Theorem \ref{cheeger's_cycles}]
Take $\epsilon >0$ sufficiently small. The set $Y\cap
S(x_0,\epsilon)$ divides $X\cap S(x_0,\epsilon)$ by two connected
components $V_1$ and $V_2$. Let $\xi$ be the chain constructed by
union of $Y\cap B[x_0,\epsilon]$ and $V_1$. Since
$\xi=\partial(W_1\cap B[x_0,\epsilon])$ we obtain that $\xi$ is a
cycle. Let us prove that $\xi$ defines a nontrivial element in
$MH_{loc,m-1}^{\nu}(X,0)$. Observe that $\xi$ is admissible because
\begin{eqnarray*}
\mu(supp(\xi),x_0)&=&\mu(supp(Y),x_0) \\
&>&\nu .
\end{eqnarray*}
If $\xi=\partial\eta$ for some chain $\eta$, then
$$W_1\cap B[x_0,\epsilon]\subset supp(\eta)\subset X$$
and
$$\mu(W_1,x_0)\leqslant \mu(supp(\eta),x_0)\leqslant \mu(X,x_0).$$
Since
$$\mu(X,x_0)=\mu(W_1,x_0)=k,$$ we get $\mu(supp(\eta),x_0)=k$,
i.e. $\eta$ is not an admissible chain. Thus, we conclude that
$[\xi]\neq 0$ in $MH_{loc,m-1}^{\nu}(X,0)$.

Now, let us prove that $c\xi$; $c\in\R$ is not homologous to the
element of $MH_{loc,m-1}^{\nu}(X,0)$ defined by the fundamental
cycle of $X\cap S(x_0,\epsilon)$. We have that $supp(f-c\xi)$ is
union of $Y\cap B[x_0,\epsilon]$ and $V_2$. Let
$f-c\xi=\partial\eta$ for some chain $\eta$. Then
$$W_2\cap B[x_0,\epsilon]\subset supp(\eta)\subset X$$
and
$$\mu(W_2,x_0)\leqslant \mu(supp(\eta),x_0)\leqslant \mu(X,x_0).$$
Since,
$$\mu(X,x_0)=\mu(W_2,x_0)=k,$$ we get $\mu(supp(\eta),x_0)=k$,
i.e. $\eta$ is not an admissible chain.
\end{proof}

The cycle $\xi$ is called \emph{Cheeger's cycle} and the set $Y$ is
called \emph{the base of the Cheeger's cycle.}

\begin{ex} Let $X\subset\R^4$ be defined by $(x_1,x_2,x_3,t)
\in\R^4$ such that
$$((x_1-t)^2+x_2^2+x_3^2-t^2)((x_1+t)^2+x_2^2+x_3^2-t^2)=t^{4\beta};
\ t\geqslant 0$$ where $\beta>2$ is a rational number. Let $Y$ be
defined by $(x_1,x_2,x_3,t)\in X$ such that $x_1=0.$ Then $Y\subset
X\subset\R^4$ are semialgebraic subsets such that
\begin{enumerate}
\item The link of $X$ at $0$ is homeomorphic to $2$-sphere $S^2$;
\item $0\in Y$ and $X-Y$ has tow connected components $W_1$ and
$W_2$;
\item $\mu(X,0)=\mu(W_1,0)=\mu(W_2,0)=3$;
\item $\mu(Y,0)=\beta + 1$.
\end{enumerate}
From above theorem, for $3<\nu <\beta +1$, $MH_{loc,2}^\nu(X,x_0)$
contains a subspace isomorphic to $\R^2$. Hence, the filtration
theorem of \cite{BB2} is not valid for $k>1$.
\end{ex}

\section{Main Results}

We say a subset $X\subset\C^3$ is a surface defined by a real
polynomial $f(x,y,z)$ when $f(x,y,z)$ is a polynomial with real
coefficients and $$X=\{(x,y,z)\in\C^3 \ : \ f(x,y,z)=0\}.$$ In this
case, we define $$X(\R)=\{(x,y,z)\in\R^3 \ : \ f(x,y,z)=0\}.$$

In the following, let $X\subset\C^3$ be a surface defined by a real
polynomial $f(x,y,z)$ which is weight homogeneous with respect to
$w=(w_1,w_2,w_3),$ where $w_1,w_2,w_3$ are coprime positive
integers. Let us suppose that $0$ is an isolated singular point of
$X$. Let $(\varphi,X^{\prime})$ be the homogeneous modification of
$X$ and $(\pi,M)$ be the Seifert Fibration associated to
$w=(w_1,w_2,w_3)$. Let $\phi\colon M^{\prime}\rightarrow M$ be the
branched covering such that
$$\pi\circ\varphi(z_1,z_2,z_3)=\phi[z_1:z_2:z_3]$$
for all $(z_1,z_2,z_3)\in X^{\prime}-\{0\}$ and let
$$F=\{[z_1:z_2:z_3]\in M^{\prime} \ : \ z_1z_2z_3=0\}.$$

\begin{prop}\label{prop_B}
Let $M^{\prime}-\phi(F)$ be with the holomorphic structure presented
in Proposition \ref{prop_A}. Then
\begin{enumerate}
\item there exists an antiholomorphic involution
$$j\colon M-\phi(F)\rightarrow M-\phi(F)$$
such that
$$\pi(X(\R)-\{0\})=\{m\in M-\phi(F) \ : \ j(m)=m\};$$
\item if $M-\pi(X(\R)-\{0\})$ is not connected, then
it has two connected components $M_1$ and $M_2$ such that
$M_1=j(M_2).$
\end{enumerate}
\end{prop}

\begin{proof} Since $X$ is defined by a real polynomial, then the
complex conjugation $$\tau\colon\C^3\rightarrow\C^3,$$ given by
$\tau(x_1,x_2,x_3)=(\overline{x_1},\overline{x_2},\overline{x_3})$,
defines a complex involution $$J=\tau_{| X}\colon\rightarrow X$$ on
$X$. We define a map $$ j\colon M\rightarrow M $$ in the following
way: given $m\in M$, let $x$ be a point on $X-\{0\}$ such that
$\pi(x)=m$, then $j(m):=\pi(J(x)).$ It is clear that
$$j\colon M\rightarrow M $$ is well defined and $j$ is an
involution. Using the holomorphic coordinates defined in Proposition
\ref{prop_A} one can show that $j$ is antiholomorphic on
$M-\phi(F)$. Now, let us show that
$$\pi(X(\R)-\{0\})=\{m\in M-\phi(F) \ : \ j(m)=m\}.$$
It is clear that $$\pi(X(\R)-\{0\})\subset\{m\in M-\phi(F) \ : \
j(m)=m\}.$$ So, let $m\in M-\phi(F)$ be such that $j(m)=m$, i.e.
$m=\pi(x_1,x_2,x_3);$ $(x_1,x_2,x_3)\in X$ and $x_1x_2x_3\neq 0.$
Since $j(m)=m$, there exists a $t\in\C^*$ such that
$$(\overline{x_1},\overline{x_2},\overline{x_3})
=(t^{w_1}x_1,t^{w_2}x_2,t^{w_3}x_3).$$ Let $s\in\C^*$ be such that
$$\overline{s}=ts.$$ Then, for each $k=1,2,3$, we have
\begin{eqnarray*}
\overline{s^{-w_k}x_k}&=& (\overline{s})^{-w_k}\overline{x_k} \\
&=& (ts)^{-w_k}\overline{x_k} \\
&=& s^{-w_k}x_k
\end{eqnarray*}
i.e. $s^{-w_k}x_k\in\R$ for all $k=1,2,3$ and
$$\alpha_w(s,(s^{-w_1}x_1,s^{-w_2}x_2,s^{-w_3}x_3)=(x_1,x_2,x_3),$$
i.e. $m\in\pi(X(\R)-\{0\}).$

Finally, let us suppose that $M-\pi(X(\R)-\{0\})$ is not connected.
Since $$\pi(X(\R)-\{0\})=\{m\in M-\phi(F) \ : \ j(m)=m\},$$ we can
compile  the proof of Proposition 5.4.2 of \cite{BR} (page 260) to
show that $M-\pi(X(\R)-\{0\})$ has two connected components $M_1$
and $M_2$ such that $M_1=j(M_2).$
\end{proof}

\begin{cor}\label{corollary} If $M-\pi(X(\R)-\{0\})$ is not connected, then
$X-{\pi}^{-1}(\pi(X(\R)))$ has exactly two connected components
$X_1,X_2$ such that $X_1=\tau(X_2)$, where
$\tau\colon\C^3\rightarrow\C^3$ is the complex conjugation
$\tau(x_1,x_2,x_3)=(\overline{x_1},\overline{x_2},\overline{x_3}).$
\end{cor}

\begin{teor}\label{theor1} Let $X\subset\C^3$ be a irreducible surface
defined by a real weighted homogeneous polynomial $f(x_1,x_2,x_3);$
with respect to $w=(w_1,w_2,w_3)$, where $w_1,w_2,w_3$ are coprime
positive integers. Suppose that the singularity of $X$ at $0\in\C^3$
is isolated. If
\begin{enumerate}
\item $2w_3 < \inf\{w_1,w_2\}$;
\item $X(\R)\neq\{0\}$;
\item $X(\R)\cap\{(x_1,x_2,x_3)\in\R^3 \ : \ x_1x_2x_3=0\}=\{0\}$.
\end{enumerate}
Then the singularity of $X$ at $0\in\C^3$ is not strongly metrically
conic.
\end{teor}

\begin{proof}
Let $i\colon X(\R)-\{0\}\rightarrow X-\{0\}$ be the embedding
induced by inclusion $X(\R)\subset X$.

\bigskip
\noindent{\tt Case 1.} We suppose that there exists a connected
component $C$ of $X(\R)-\{0\}$ such that $i(C)$ presents nontrivial
element in $H_1(X-\{0\})$. Let $Y=\overline{C}=C\cup\{0\}$. By the
results of \cite{BF}, $Y$ is bi-Lipschitz equivalent, with respect
to inner metric, to a $\beta$-horn $H_{\beta}$, where
$\displaystyle\beta =\frac{\inf\{w_1,w_2\}}{w_3}$ and from
hypothesis $\beta>1$. We found a 1-dimensional nontrivial cycle
$\sigma$ in $H_1(X-\{0\})$, given by $\sigma =\partial\eta$ where
$\eta =Y\cap B[0,\epsilon]$. It means that 1-dimensional
characteristic exponent of $X$ at $0$ (see \cite{BC}, \cite{BB1},
\cite{BB2}) is bigger than or equal to $\mu(\eta,0)=\beta+1>2$.
Therefore, by results of \cite{BB2}, $X$ at $0$ is not strongly
metrically conic, because otherwise the corresponding exponent must
be smaller than or equal to $2$.

\bigskip
\noindent{\tt Case 2.} We suppose $i(X(\R)-\{0\})$ presents
nontrivial element in the homology group $H_1(X-\{0\})$. Let
$\gamma=\pi(X(\R)-\{0\})$, where $(\pi,M)$ is the Seifert Fibration
of $X-\{0\}$ associated to $w=(w_1,w_2,w_3)$. Then, $[\gamma]$ is a
trivial element in $H_1(M).$ Since $M$ is a $2$-dimensional
orbifold, we obtain that $M-\gamma$ is not connected. Using
Proposition \ref{prop_B}, we obtain that $M-\gamma$ has two
connected components $M_1$ and $M_2$ such that $M_1=j(M_2).$
Moreover, by Corollary \ref{corollary}, we obtain that
$X-{\pi}^{-1}(\pi(X(\R)))$ is a union of two connected components
$X_1$ and $X_2$ such that $X_1=\tau(X_2)$, where $\tau$ is the
complex conjugation of $\C^3$. The set $Y=X-{\pi}^{-1}(\pi(X(\R)))$
is obtained by the revolution of $X(\R)$ by a $1$-dimensional
subgroup of isometry group of $\C^3$. Since $\mu(X(\R),0)=\beta +1$
, where $\displaystyle\beta=\frac{\inf\{w_1,w_2\}}{w_3}$ (see
\cite{BF}), we obtain that
\begin{eqnarray*}
\mu(Y,0)&=&\mu(X(\R),0)+1 \\
&=&\beta +2
\end{eqnarray*}
and, since $\beta>2$, we have $\mu(Y,0)>4.$ From the other hand,
$$\mu(X,0)=\mu(X_1,0)=\mu(X_2,0)=4.$$
Now, since $X-\{0\}$ is connected one can apply Theorem
\ref{cheeger's_cycles} and, by the remark, $X$ is not strongly
metrically conic at the singular point $0$.
\end{proof}

\begin{teor}\label{theor2} Let $X\subset\C^3$ be a irreducible surface
defined by a real weighted homogeneous polynomial $f(x_1,x_2,x_3);$
with respect to $w=(w_1,w_2,w_3),$ where $w_1,w_2,w_3$ are coprime
positive integers. Suppose that the singularity of $X$ at $0\in\C^3$
is isolated. If
\begin{enumerate}
\item $w_3<\inf\{w_1,w_2\}$;
\item $\pi(X(\R)\subset M$ has more than one
connected component;
\item $X(\R)\cap\{(x_1,x_2,x_3)\in\R^3 \ : \ x_1x_2x_3=0\}=\{0\}$.
\end{enumerate}
Then the singularity of $X$ at $0\in\C^3$ is not metrically conic.
\end{teor}

In order to show this theorem, we need the following proposition.

\begin{prop} Let $Y$ be a connected component of $X(\R)-\{0\}$ and
let $\xi =Y\cap S(0,\epsilon)$, for sufficiently small $\epsilon>0$.
Let $[\xi]\neq 0$ in $H_1(X-\{0\})$. Then the singularity of $X$ at
$0$ is not metrically conic.
\end{prop}

\begin{proof} Suppose that there exist a subset
$$N\subset\{x\in\R^m \ : \ \| x\|=1\}$$ and a bi-Lipschitz map-germ
$$F\colon (X,0)\rightarrow (C_0N,0).$$

Given $r>0$, sufficiently small, let $\xi_r=Y\cap S(0,r).$ By the
conditions of the theorem, $[\xi_r]\neq 0$ in $H_1(X-\{0\}).$ Hence
$F_{*}[(\xi_r)]\neq 0$ in $H_1(C_0N-\{0\})$. Let us denote by
$Band(k_1,k_2)$ the following subset
$$Band(k_1,k_2)=\{x\in\R^m \ : \ k_1\leqslant \|x\| \leqslant k_2\},$$
where $k_1<k_2$ are bi-Lipschitz constants of $F$.  Let
$$\theta_r\colon\R^m-\{0\}\rightarrow\R^m-\{0\}$$ be defined by
$$\theta_r(x)=\frac{1}{r}x.$$
Since $F$ is a bi-Lipschitz map-germ with bi-Lipschitz constants
$k_1<k_2$, we obtain $$\theta_r(F(\xi_r))\subset Band(k_1,k_2)\cap
C_0N$$ and
$$diam(\theta_r(F(\xi_r)))\leqslant\frac{1}{r}k_2diam(\xi_r).$$ On
the other hand, since the germ $(Y,0)$ is bi-Lipschitz equivalent to
a $\beta$-horn, where $\beta=\frac{\inf\{w_1,w_2\}}{w_3}$ (see
\cite{BF}), we obtain
$$diam(\xi_r)\leqslant \tilde{k}r^{\beta}$$ for some constant
$\tilde{k}>0$. Thus,
$$diam(\theta_r(F(\xi_r)))\leqslant k_2\tilde{k}r^{\beta-1}$$
and, in particular, it means that
$$\lim_{r\to 0}diam(\theta_r(F(\xi_r)))=0.$$

Let $P\colon Band(k_1,k_2)\cap C_0N\rightarrow N$ be a canonical
projection
$$P(x)=\frac{1}{\|x\|}x.$$ Since $P$ is a Lipschitz map, we have
$$\lim_{r\to 0}P(\theta_r(F(\xi_r)))=0.$$ Since $N$ is a topological
manifold, $P(\theta_r(F(\xi_r)))$ defines a trivial cycle in
$H_1(N)$, for sufficiently small $r>0$. This is a contradiction
because $P_*$ and $(\theta_r)_*$ are isomorphisms.
\end{proof}

\begin{proof}[Proof of the theorem.]
Suppose that $\pi(X(\R))$ is trivial on $M.$ Since $M-\pi(X(\R))$
has more than one connected component, there exists a component
$C\subset\pi(X(\R))$ such that $[C]\neq 0$ in $H_1(M).$ Then
$Y={\pi}^{-1}(C)$ satisfies the conditions of proposition above.
This proves the theorem.
\end{proof}

\section{Example: Surfaces of Brieskorn}

\begin{teor} The singularity at $0\in\C^3$ of Brieskorn surface $X$
defined as follows:
$$x^2+y^2=z^{2k}; \quad k>2$$ is not strongly metrically conic.
\end{teor}

Note that these surfaces do not satisfy the conditions of Theorem
\ref{theor1} and the conditions of Theorem \ref{theor2}.

\begin{proof}{Proof of the theorem.} Let $Y\subset X$ be the result
of the $\C^{*}$-action on $X(\R)$. Let us show that $Y$ is a base of
a Cheeger's cycle on $X$. Let us show that $X-Y$ has exactly two
connected components. Consider the hyperplane section
$\hat{X}=X\cap\{z=1\}$. It is an affine curve given by the equation
$$ x^2+y^2=1.$$ The set $\hat{X}(\R)$ is homeomorphic to $S^1$ and
$\hat{X}$ is homeomorphic to a cone over $S^1$, hence the set
$\hat{X}-\hat{X}(\R)$ contains exactly two  connected components,
which are conjugated. Let $(x_0,y_0,z_0)\in X-Y$. We are going to
show that there is no continuous path connecting $(x_0,y_0,z_0)$
with $(\overline{x_0},\overline{y_0},\overline{z_0})$. Observe that
we can suppose that $(x_0,y_0,z_0)$ does not belong to
$X\cap\{z=0\}$, otherwise we can take a nearly point belonging to
$X-\{z=0\}$ connected to this one by a continuous path on $X-Y$)
Suppose that there exists a continuous path $\gamma\colon
[0,1]\rightarrow X$ such that $\gamma(0)=p$ and
$\gamma(1)=\overline{p}$. By a transversality argument, we can
suppose that $\gamma$ does not intersect the set $X\cap\{z=0\}$. Let
$$\rho\colon X-\{z=0\}\rightarrow\hat{X}$$ defined by
$\rho(x,y,z)=(z^{-k}x,z^{-k}y,1)$. Since the map $\rho$ respects the
complex conjugation, we obtain $$\overline{\rho(x_0,y_0,z_0)}=
\rho(\overline{x_0},\overline{y_0},\overline{z_0}).$$ Thus
$\rho(x_0,y_0,z_0)$ and
$\rho(\overline{x_0},\overline{y_0},\overline{z_0})$ belong to the
different connected components of $\hat{X}-\hat{X}(\R)$. Hence the
path $\gamma$ must intersect $Y$, because
$\rho(Y-{z=0})=\hat{X}(\R)$. It means that $X-Y$ is not connected.
Moreover $X-Y$ has exactly two connected components $X_1$ and $X_2$
which are conjugated. Now, since $\mu(X(\R),0)=k+1$, we get that
$\mu(Y,0)\geqslant k+2$, using the same argument as in the proof of
Theorem \ref{theor1}.

Finally, we obtain that the sets $X$ and $Y$ satisfy the conditions
of Theorem \ref{cheeger's_cycles}. The set $Y$ is a base of a
Cheeger's cycle and the singular point $\{0\}$ is not strongly
metrically conic.
\end{proof}

\end{document}